\documentclass[runningheads,a4paper]{llncs}

\usepackage{amssymb}
\usepackage{graphicx}

\usepackage{amsmath,color,subfig}

\usepackage{algorithm}
\usepackage{algorithmic}

\usepackage{url}
\urldef{\mailsa}\path|{roberto.cavoretto, alessandra.derossi, emma.perracchione}@unito.it|
\newcommand{\keywords}[1]{\par\addvspace\baselineskip
\noindent\keywordname\enspace\ignorespaces#1}

\newcommand{\NN}{\mathbb{N}}

\newcommand{\RR}{\mathbb{R}}

\begin{document}

\mainmatter  

\title{Partition of Unity Interpolation on Multivariate Convex Domains}

\titlerunning{Partition of Unity Interpolation on Multivariate Convex Domains}

%
%
\author{Roberto Cavoretto
\and Alessandra De Rossi\and Emma Perracchione}
\authorrunning{R. Cavoretto\and A. De Rossi\and E. Perracchione}

\institute{Department of Mathematics \lq\lq G. Peano\rq\rq, University of Torino,\\
Via Carlo Alberto 10, 10123 Torino, Italy\\
\mailsa\\
\url{}}

%
%

\toctitle{Lecture Notes in Computer Science}
\tocauthor{R. Cavoretto, A. De Rossi, E. Perracchione}
\maketitle

\begin{abstract}
In this paper we present a new algorithm for multivariate interpolation of scattered data sets lying in convex domains $\Omega \subseteq \RR^N$, for any $N \geq 2$. To organize the points in a multidimensional space, we build a $kd$-tree space-partitioning data structure, which is used to efficiently apply a partition of unity interpolant. This global scheme is combined with local radial basis function approximants and compactly supported weight functions. A detailed description of the algorithm for convex domains and a complexity analysis of the computational procedures are also considered. Several numerical experiments show the performances of the interpolation algorithm on various sets of Halton data points contained in $\Omega$, where $\Omega$ can be any convex domain like a 2D polygon or a 3D polyhedron.\keywords{Meshfree Approximation, Multivariate Algorithms, Partition of Unity Methods, Scattered Data.}
\end{abstract}

\section{Introduction}


In this paper we deal with the problem of interpolating a (usually) large number of multivariate scattered data points lying in convex domains or, more precisely, in convex hulls $\Omega \subseteq \RR^N$, for any $N \geq 2$. In general, this problem is considered in literature supposing to interpolate data points which are situated in suitable or simple domains such as hypercubes or hyperrectangles (see e.g. \cite{Fasshauer07,Fasshauer11,Nguyen08}). Thus we construct a numerical algorithm which can efficiently be used for scattered data interpolation in $\Omega$. To organize the points in a multivariate space, we make use of a space-partitioning data structure known as \textsl{$kd$-tree} (see \cite{Samet90}). This code is designed to numerically approximate data points by the partition of unity method, a global interpolation scheme which is combined with local radial basis function (RBF) approximants and compactly supported weight functions (see \cite{Iske11,Wendland02,Wendland05}). A detailed design of this algorithm as well as an analysis of its complexity is considered. 

Moreover, we observe that the implemented code is completely automatic and any choice depending on the space dimension has suitably been studied so that this algorithm can work for any dimension. Numerical experiments show the performances of the interpolation algorithm on various sets of Halton data points contained in $\Omega \subseteq \RR^N$, for $N = 2,3$. Here, $\Omega$ is any convex domain like a 2D polygon (e.g., a triangle or a hexagon) or a 3D polyhedron (e.g., a pyramid or a cylinder). Note that this algorithm for convex hulls extends our previous works on the topic \cite{Cavoretto12,Cavoretto14a,Cavoretto14b,Cavoretto14c}.

The paper is organized as follows. In Section \ref{PU-RBF} we give a general presentation of the partition of unity interpolation combined with local radial basis functions, reporting some theoretical results. In Section \ref{alg}, we describe the algorithm for convex hulls and analyze its complexity. In Section \ref{num_res}, in order to show accuracy and efficiency of the interpolation algorithm, we report numerical experiments considering various sets of scattered data points contained in 2D and 3D convex domains. Finally, Section \ref{concl} refers to conclusions and future work.


\section{Partition of Unity Interpolation} \label{PU-RBF}
Let ${\cal X}_n=\{\boldsymbol{x}_i, i = 1,2, \ldots , n \}$ be a set of distinct data points, arbitrarily distributed in a domain $\Omega \subseteq \RR^N$, $N \geq 1$, with an associated set ${\cal F}_n=\{f_i, i = 1,2,  \ldots , n \}$ of data values or function values, which are obtained by sampling some (unknown) function $f:\Omega \rightarrow \RR$ at the data points, i.e., $f_i=f(\boldsymbol{x}_i)$, $i=1,2,\ldots,n$. 

The basic idea of the partition of unity method is to start with a partition of the open and bounded domain $\Omega \subseteq \RR^N$ into $d$  subdomains $\Omega_j$ such that $\Omega \subseteq \bigcup_{j=1}^{d} \Omega_j$ with some mild overlap among the subdomains. Associated with these subdomains we choose a partition of unity, i.e. a family of compactly supported, non-negative, continuous functions $W_j$ with $\text{supp}(W_j) \subseteq \Omega_j$ such that 
\begin{equation}
\sum_{j=1}^{d} W_j(\boldsymbol{x}) = 1, \hspace{1cm} \boldsymbol{x} \in \Omega.
\end{equation}
The global approximant thus assumes the following form
\begin{eqnarray}
\label{pui}
	{\cal I}(\boldsymbol{x})= \sum_{j=1}^{d} R_j(\boldsymbol{x}) W_j(\boldsymbol{x}), \hspace{1cm} \boldsymbol{x} \in \Omega.
\end{eqnarray}
For each subdomain $\Omega_j$ we define a local RBF interpolant $R_j:\Omega \rightarrow \RR$ of the form 
\begin{eqnarray}
\label{intfun}
R_j({\boldsymbol{x}})=\sum_{k=1}^{n_j} c_{k} \phi (d(\boldsymbol{x},\boldsymbol{x}_{k})),
\end{eqnarray}
where $d(\boldsymbol{x},\boldsymbol{x}_k)=||\boldsymbol{x}-\boldsymbol{x}_k||_2$ is the Euclidean distance, $\phi:[0,\infty) \rightarrow \RR$ is called \textsl{radial basis function}, and $n_j$ indicates the number of data points in $\Omega_j$. Moreover, $R_j$ satisfies the interpolation conditions
\begin{eqnarray}
\label{condinterp} 
  R_j(\boldsymbol{x}_i)=f_i, \hspace{1.cm} i=1,2,\ldots,n_j.
\end{eqnarray} 
In particular, we observe that if the local approximants satisfy the interpolation conditions \eqref{condinterp}, then the global approximant also interpolates at $\boldsymbol{x}_i$, i.e. ${\cal I}(\boldsymbol{x}_i)=f(\boldsymbol{x}_i)$, for $i=1,2,\ldots,n_j$.  

Solving the $j$-th interpolation problem \eqref{condinterp} leads to a system of linear equations of the form
\begin{equation} \label{mat}
 \Phi\boldsymbol{c}=\boldsymbol{f}, \nonumber
\end{equation}
where entries of the interpolation matrix $\Phi$ are
\begin{eqnarray}
\Phi_{ik}=\phi(d(\boldsymbol{x}_i,\boldsymbol{x}_k)), \hskip 1.cm i,k=1,2,\ldots,n_j, \nonumber
\end{eqnarray}
$\boldsymbol{c}=[c_1,c_2,\ldots,c_{n_j}]^T$ and $\boldsymbol{f}=[f_1,f_2,\ldots,f_{n_j}]^T$.

Now, we give the following definition (see \cite{Wendland02}).
\begin{definition}
Let $\Omega \subseteq \RR^N$ be a bounded set. Let $\{\Omega\}_{j=1}^{d}$ be an open and bounded covering of $\Omega$. This means that all $\Omega_j$ are open and bounded and that $\Omega$ is contained in their union. A family of nonnegative functions $\{W_j\}_{j=1}^{d}$ with $W_j \in C^k(\RR^N)$ is called a $k$-stable partition of unity with respect to the covering $\{\Omega_j\}_{j=1}^{d}$ if
\begin{enumerate}
	\item[1)] ${\rm supp}(W_j) \subseteq \Omega_j$;
	\item[2)] $\sum_{j=1}^{d} W_j(\boldsymbol{x}) \equiv 1$ on $\Omega$;
	\item[3)] for every $\beta \in \NN_0^N$ with $|\beta| \leq k$ there exists a constant $C_{\beta} > 0$ such that
\begin{equation}
	||D^{\beta}W_j||_{L_{\infty}(\Omega_j)}\leq C_{\beta}/\delta_j^{|\beta|}, \hspace{1.cm} j=1,2,\ldots,d, \nonumber
\end{equation}
where $\delta_j = {\rm diam}(\Omega_j)=\sup_{\boldsymbol{x},\boldsymbol{y} \in \Omega_{j}} ||\boldsymbol{x}-\boldsymbol{y}||_2$.
\end{enumerate}
\end{definition}

In accordance with the statements in \cite{Wendland02} we require some additional regularity assumptions on the \textsl{covering} $\{\Omega_j\}_{j=1}^{d}$.


\begin{definition} \label{defpr}
Suppose that $\Omega \subseteq  \RR^N$ is bounded and ${\cal X}_n=\left\{\boldsymbol{x}_i, i=1,2,\ldots,n\right\} \subseteq \Omega$ are given. An open and bounded covering $\{\Omega_j\}_{j=1}^{d}$ is called regular for $(\Omega,{\cal X}_n)$ if the following properties are satisfied:
\begin{itemize}
	\item[(a)] for each $\boldsymbol{x} \in \Omega$, the number of subdomains $\Omega_j$ with $\boldsymbol{x} \in \Omega_j$ is bounded by a global constant $K$;
	\item[(b)] each subdomain $\Omega_j$ satisfies an interior cone condition \cite{Wendland05};
	\item[(c)] the local fill distances $h_{{\cal X}_{n_j}, \Omega_j}$, where ${\cal X}_{n_j}={\cal X}_n \cap \Omega_j$, are uniformly bounded by the global fill distance $h_{{\cal X}_n, \Omega}$, i.e.	
	\begin{eqnarray} 
		h_{{\cal X}_{n}, \Omega} = \sup_{\boldsymbol{x} \in \Omega}\min_{\boldsymbol{x}_k\in {\cal X}_n} d(\boldsymbol{x},\boldsymbol{x}_k). \nonumber
	\end{eqnarray}	
\end{itemize}
\end{definition}

After defining the space $C_{\nu}^k(\RR^N)$ of all functions $f \in C^k$ whose derivatives of order $|\beta|=k$ satisfy $D^{\beta}f(\boldsymbol{x})= {\cal O}(||\boldsymbol{x}||_2^{\nu})$ for $||\boldsymbol{x}||_2 \rightarrow 0$, we consider the following convergence result (see, e.g., \cite{Fasshauer07,Wendland05}).

\begin{theorem}
	Let $\Omega \subseteq  \RR^N$ be open and bounded and suppose that ${\cal X}_n = \{\boldsymbol{x}_i, i=1,$ $2,\ldots,n \}\subseteq \Omega$. Let $\phi \in C_{\nu}^k(\RR^N)$ be a strictly positive definite function. Let $\{\Omega_j\}_{j=1}^{d}$ be a regular covering for $(\Omega, {\cal X}_n)$ and let $\{W_j\}_{j=1}^{d}$ be $k$-stable for $\{\Omega_j\}_{j=1}^{d}$. Then the error between $f \in {\cal N}_{\phi}(\Omega)$, where ${\cal N}_{\phi}$ is the native space of $\phi$, and its partition of unity interpolant (\ref{pui}) can be bounded by
\begin{equation}
	|D^{\beta}f(\boldsymbol{x}) - D^{\beta}{\cal I}(\boldsymbol{x})| \leq C h_{{\cal X}_n, \Omega}^{(k+\nu)/2 - |\beta|} |f|_{{\cal N}_{\phi}(\Omega)}, \nonumber
\end{equation}
for all $\boldsymbol{x} \in \Omega$ and all $|\beta| \leq k/2$. 
\end{theorem}  

If we compare this result with the global error estimates (see e.g. \cite{Wendland05}), we can see that the partition of unity preserves the local approximation order for the global fit. This means that we can efficiently compute large RBF interpolants by solving small RBF interpolation problems and then glue them together with the global partition of unity $\{W_j\}_{j=1}^{d}$. In other words, the partition of unity approach is a simple and effective technique to decompose a large problem into many small problems while at the same time ensuring that the accuracy obtained for the local fits is carried over to the global one. 



\section{Algorithm for Convex Hulls} \label{alg}

In this section we present an algorithm for multivariate interpolation of scattered data sets lying in a convex domain (or convex hull) $\Omega \subseteq \RR^N$, for any $N \geq 2$. This code is based on a global partition of unity interpolant using local RBF interpolants and compactly supported weight functions. To organize the points in a multivariate space, we build an efficient space-partitioning data structure as the \textsl{kd-trees}, because this enables us to efficiently answer a query, known as \textsl{range search} (see \cite{Arya98,deBerg97}). In fact, we need to solve the following computational issue:
\begin{itemize}
	\item[] \emph{Given a set $X$ of points $\boldsymbol{x}_i \in X$ and a subdomain $\Omega_j$, find all points situated in that subdomain, i.e. $\boldsymbol{x}_i \in X_j = X \cap \Omega_j$.} 
\end{itemize} 
Note that the subdomain $\Omega_j$ denotes a generic region, so the index $j$ is here fixed. For simplicity, all details of this algorithm concern a generic convex hull $\Omega \subseteq [0,1]^N$, but its generalization is obviously possible and straightforward.

\subsection{Description of the Algorithm} \label{des_alg}

\noindent INPUT: $N$, space dimension; $n$, number of data; ${\cal X}_n=\{\boldsymbol{x}_i, i=1,2,\ldots,n\}$, set of data points; ${\cal F}_n=\{f_i, i=1,2,\ldots,n\}$, set of data values.
\vskip 4pt
\noindent OUTPUT: ${\cal A}_s=\{{\cal I}(\tilde{\boldsymbol{x}}_i), i=1,2,\ldots,s\}$, set of approximated values.
\vskip 4pt
\noindent {\tt Stage 1.} The set ${\cal X}_n$ of data points and the set ${\cal F}_n$ of data values are loaded.
\vskip 4pt
\noindent {\tt Stage 2.} After computing the number $d$ of subdomain points, a set ${\cal C}_d=\{\bar{\boldsymbol{x}}_j, j=1,2,\ldots,d\} \subseteq \Omega$ of subdomain points is constructed\footnote{This set is obtained by generating a grid of equally spaced points on the hypercube $[0,1]^N$. They are then automatically reduced taking only those in $\Omega$ by the \texttt{inhull} \textsc{Matlab} function. Such points are the centres of partition of unity subdomains.}. Note that the number $d$ depends on both the data point number $n$ and the space dimension $N$; furthermore, it is suitably chosen assuming that the ratio $n/d \approx 2^{N+1}$.
\vskip 4pt 
\noindent {\tt Stage 3.} The number $s$ of evaluation points is computed and a set ${\cal E}_s=\{\tilde{\boldsymbol{x}}_i, i=1,2,\ldots,s\} \subseteq \Omega$ of evaluation points is generated.
\vskip 4pt 
\noindent {\tt Stage 4.} For each subdomain point $\bar{\boldsymbol{x}}_j$, $j=1,2,\ldots,d$, a local spherical subdomain is constructed, whose radius is
\begin{eqnarray}
\label{delta}
\delta_{\Omega_j} = \frac{\sqrt{2}}{D^{1/N}},
\end{eqnarray}
where $D$ is the number of subdomain points initially generated on $[0,1]^N$.
\vskip 4pt
\noindent {\tt Stage 5.} The $kd$-tree data structures are built for the set ${\cal X}_n$ of data points and the set ${\cal E}_n$ of evaluation points. 
\vskip 4pt 
\noindent {\tt Stage 6.} For each subdomain $\Omega_j$, $j=1,2,\ldots,d$, the range query problem is considered, adopting the related searching procedure which consists of the following two steps: 
\begin{enumerate}
\item[i)] Find all data points (i.e. the set ${\cal X}_{n_j}$) belonging to the subdomain $\Omega_j$ and construct a local interpolation RBF matrix by ${\cal X}_{n_j}$, where $n_j$ denotes the point number of ${\cal X}_{n_j}$.
\item[ii)] Determine all evaluation points (i.e. the set ${\cal E}_{s_j}$) belonging to the subdomain $\Omega_j$ and build a local evaluation RBF matrix by ${\cal E}_{s_j}$, where $s_j$ is the point number of ${\cal E}_{s_j}$.
\end{enumerate}
\vskip 4pt 
\noindent {\tt Stage 7.} A local RBF interpolant $R_j$ and a weight function $W_j$, $j=1,2,\ldots,d$, is computed for each evaluation point.
\vskip 4pt 
\noindent {\tt Stage 8.} The global fit (\ref{pui}) is applied, accumulating all the $R_j$ and $W_j$.
\vskip 4pt

In this algorithm for convex domains the local interpolants are computed by using compactly supported RBFs as the Wendland functions. However, this approach is completely automatic and turns out to be very flexible, since different choices of local approximants, either globally or compactly supported, are allowed.


\subsection{Complexity Analysis} \label{comp_anal}

The algorithm is based on the construction of $kd$-tree data structures. They enable us to efficiently determine all data points belonging to each subdomain $\Omega_j$, $j=1,2,\ldots,d$, so that we can compute local RBF interpolants to be used in the partition of unity scheme. Then, assuming that the covering $\{\Omega_j\}_{j=1}^{d}$ is regular and local and the set ${\cal X}_n$ of data points is quasi-uniform, we analyze the complexity of this interpolation algorithm. 

In \texttt{Stages 1-4} we have a sort of preprocessing phase where we automatically load all data sets and define the parameters concerning data, subdomain and evaluation points. To construct an algorithm which efficiently works in a generic space dimension $N$, we require that the subdomain number $d$ is proportional to the data point number $n$, taking $n/d \approx 2^{N+1}$.

In \texttt{Stage 5} we build the $kd$-trees, which needs ${\cal O}(N n\log n)$ time and ${\cal O}(N n)$ space for $n$ data points and ${\cal O}(N s\log s)$ time and ${\cal O}(N s)$ space for $s$ evaluation points. Then, in \texttt{Stage 6} we make use of the range search procedure for each subdomain $\Omega_j$, $j=1,2,\ldots,d$, whose running times are ${\cal O}(\log n)$ and ${\cal O}(\log s)$, respectively (see \cite{Wendland05}). 

Since the number of centres in each subdomain $\Omega_j$ is bounded by a constant (see Definition \ref{defpr}), we need ${\cal O}(1)$ space and time for each subdomain to solve the local RBF interpolation problems. In fact, in order to obtain the local RBF interpolants, we have to solve $d$ linear systems of (relatively) small sizes, i.e. $n_j \times n_j$, with $n_j << n$, thus requiring a constant running time ${\cal O}(n_j^3)$, $j=1,2,\ldots,d$, for each subdomain. Besides reporting the points in each subdomain in ${\cal O}(1)$, as the number $d$ of subdomains $\Omega_j$ is bounded by ${\cal O}(n)$, this leads to ${\cal O}(n)$ space and time for solving all of them. 

Thus, in \texttt{Stage 7} and \texttt{8} we have to add up a constant number of local RBF interpolants to get the value of the global fit \eqref{pui}. This can be computed in ${\cal O}(1)$ time.


\section{Numerical Experiments} \label{num_res}
In this section we present some numerical experiments we made to test our procedures implemented in \textsc{Matlab} environment. All the tests have been carried out on a Intel Core i7-4500U 1.8 GHz processor. In our results we report errors and CPU times obtained by running the algorithm on a few scattered data sets, which are located in a convex hull $\Omega \subseteq [0,1]^N$, for $N=2,3$. As interpolation points, we take uniformly random Halton data points. They are generated by using the program \texttt{haltonseq.m}, available at \cite{MCFE}, and then suitably reduced to $\Omega$. We observe that this code for convex domains is completely automatic and, though we here focus only on bivariate and trivariate interpolation, it might also be used in higher dimensions.

In order to point out accuracy of this algorithm, we compute on a reduced grid of $s$ evaluation points\footnote{The number $s$ depends on the convex domain $\Omega$; at first, we construct a uniform grid of $40^N$ points, and then we automatically reduce them taking only those in $\Omega$ through the \texttt{inhull} \textsc{Matlab} function.} Maximum Absolute Error (MAE) and Root Mean Square Error (RMSE), whose formulas are
\begin{eqnarray} \label{MAE}
	MAE = \max_{1\leq i \leq s} |f(\tilde{\boldsymbol{x}}_i) - {\cal I}(\tilde{\boldsymbol{x}}_i)|,
\end{eqnarray}
and
\begin{eqnarray} \label{RMSE}
	RMSE = \sqrt{\frac{1}{s}\sum_{i=1}^{s} |f(\tilde{\boldsymbol{x}}_i) - {\cal I}(\tilde{\boldsymbol{x}}_i)|^2}.
\end{eqnarray}
 
Moreover, we report results obtained by using as basis the Wendland $C^2$ function, i.e.,  
\begin{align*}
\phi(r) = \left(1-\delta r\right)_+^4\left(4\delta r+1\right),
\end{align*}
where $\delta \in \RR^+$ is a \textsl{shape parameter}, $r=||\cdot||_2$ is the Euclidean distance, and $(\cdot)_+$ denotes the truncated power function. We remark that this RBF is compactly supported (i.e., its support is $\left[0,1/\delta\right]$) and strictly positive definite in $\RR^N$ for $N\leq 3$ (see \cite{Wendland05}). Note that here it is used as both a basis function and a localizing function of Shepard's weight $W_j$ in the global fit \eqref{pui}.  



\subsection{Results for 2D Convex Hulls} \label{num_exp_2D}

\begin{figure}
\centering
\includegraphics[height=3.4cm]{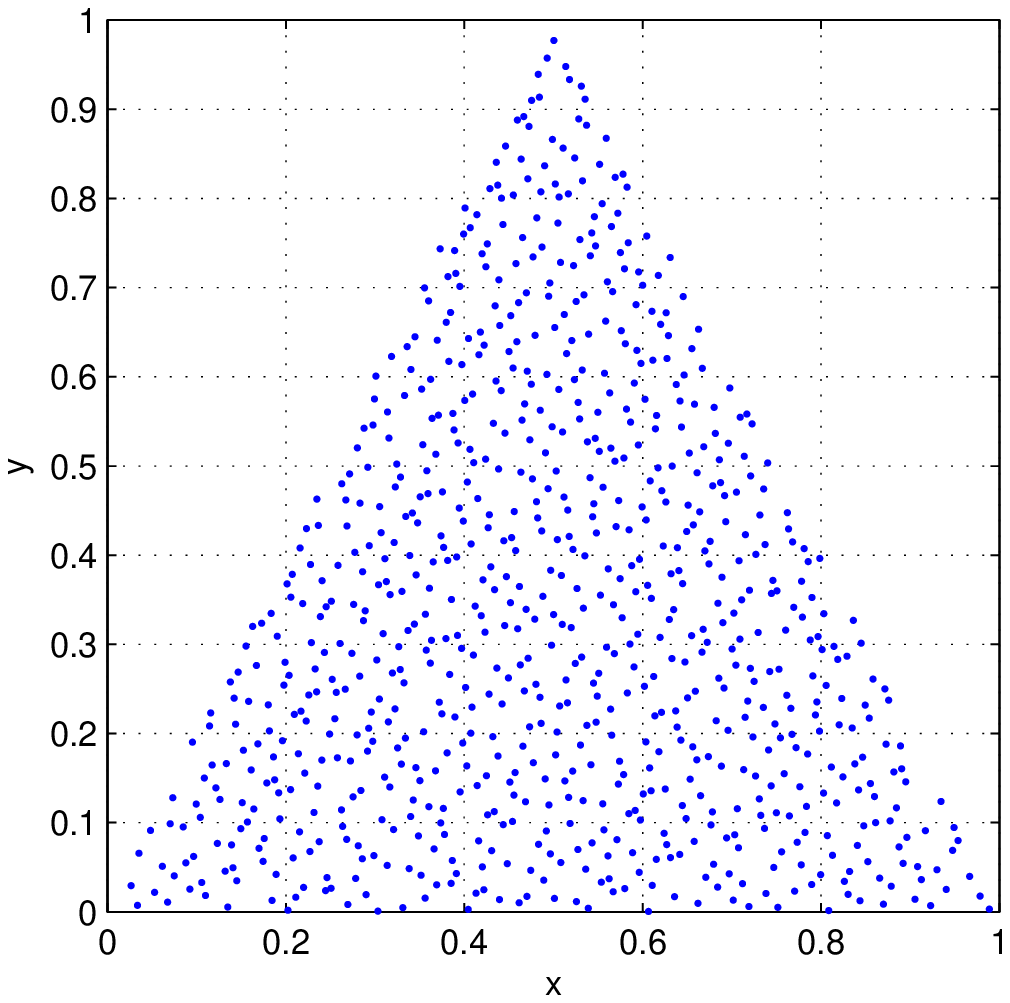} \hskip -.85cm
\includegraphics[height=3.4cm]{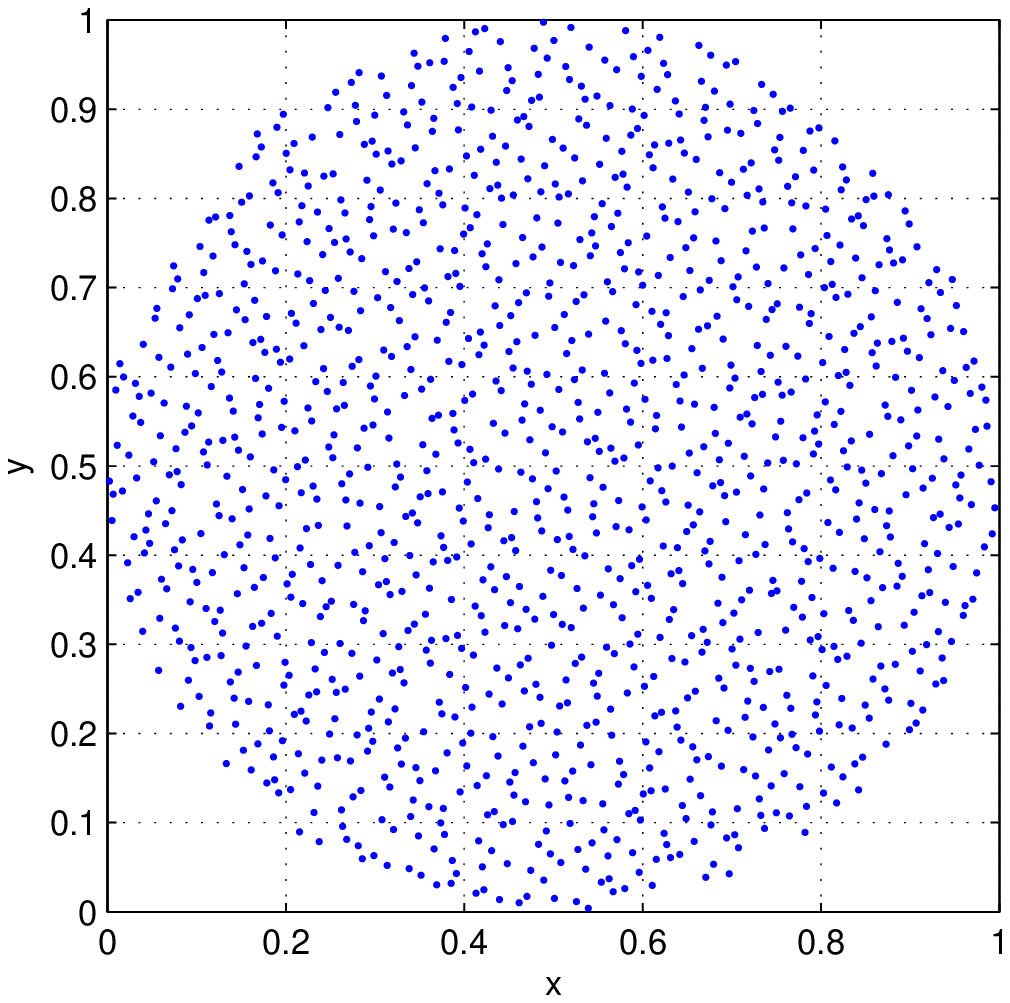} \hskip -.85cm 
\includegraphics[height=3.4cm]{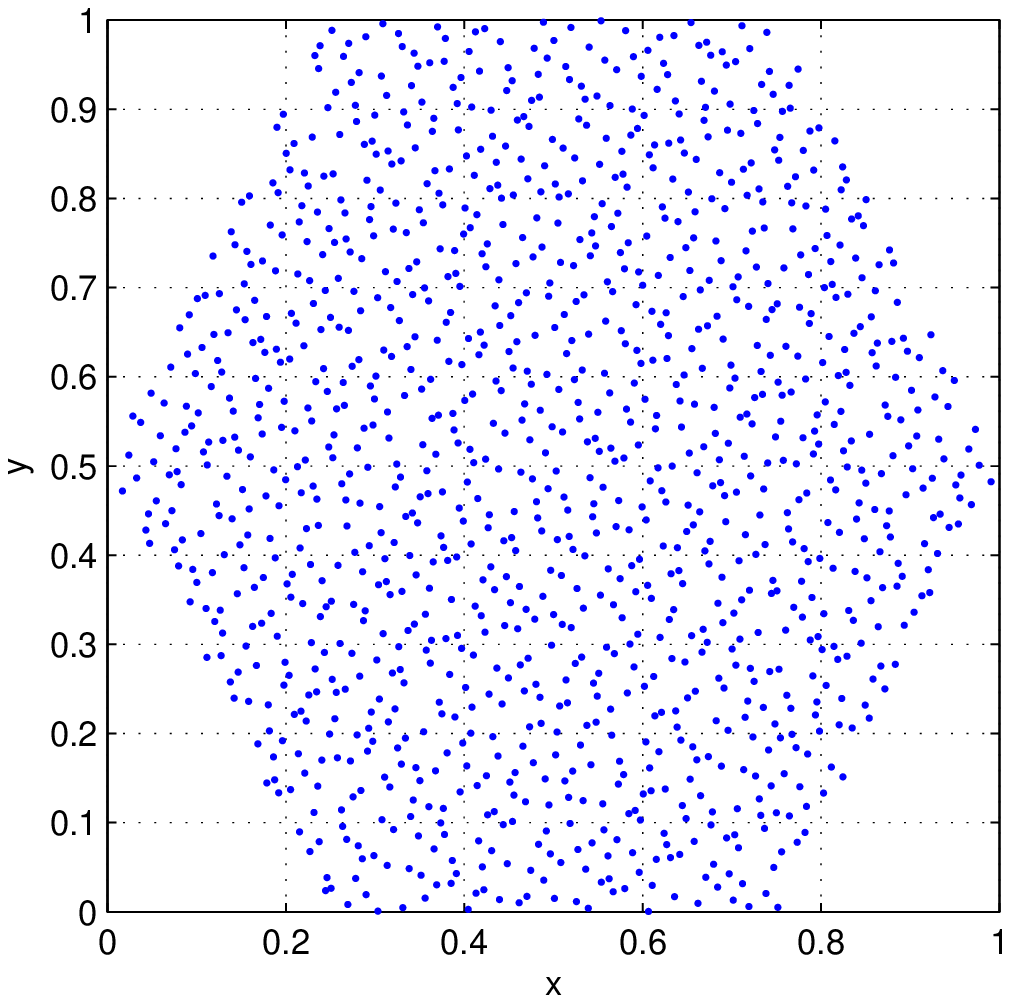}
\caption{Examples of points in 2D convex hulls. Left: triangle, 805 nodes; center: disk, 1257 nodes; right: hexagon, 1204 nodes.}
\label{fig_1}
\end{figure}

In this subsection we focus on bivariate interpolation, analyzing performances of our algorithm for convex hulls and showing the numerical results obtained by considering five sets of Halton data points. These tests are carried out considering different convex domains, i.e., a triangle, a disk and a hexagon (see Figure \ref{fig_1}).

In the various experiments we investigate accuracy of the interpolation algorithm taking the data values by the well-known 2D Franke's test function
\begin{align*}
f_2(x_1,x_2)&= \frac{3}{4}{\rm e}^{-\frac{(9x_1-2)^2+(9x_2-2)^2}{4}}+\frac{3}{4}{\rm e}^{-\frac{(9x_1+1)^2}{49}-\frac{9x_2+1}{10}} \\
&+\frac{1}{2} {\rm e}^{-\frac{(9x_1-7)^2+(9x_2-3)^2}{4}}-\frac{1}{5} {\rm e}^{-(9x_1-4)^2-(9x_2-7)^2}.
\end{align*}

\begin{figure}
\centering
\includegraphics[height=3.4cm]{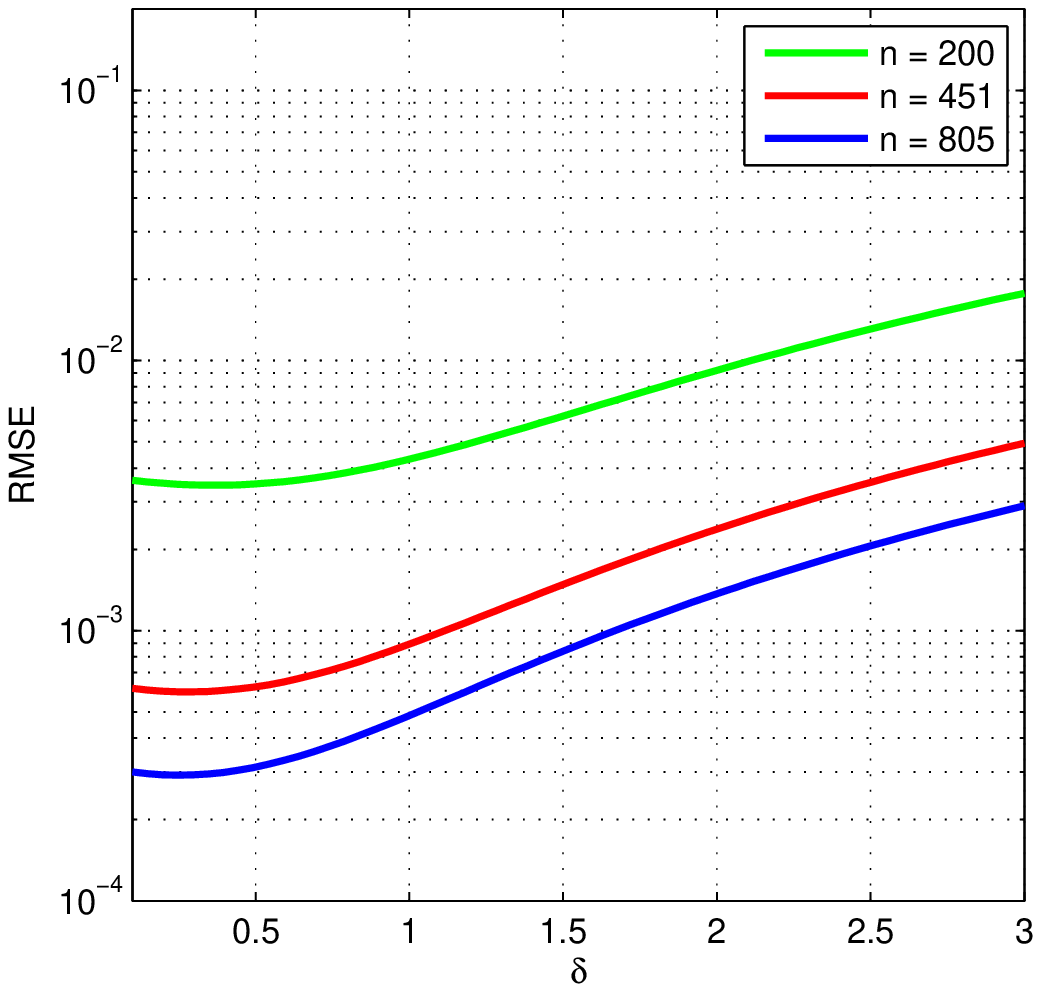} \hskip -0.85cm
\includegraphics[height=3.4cm]{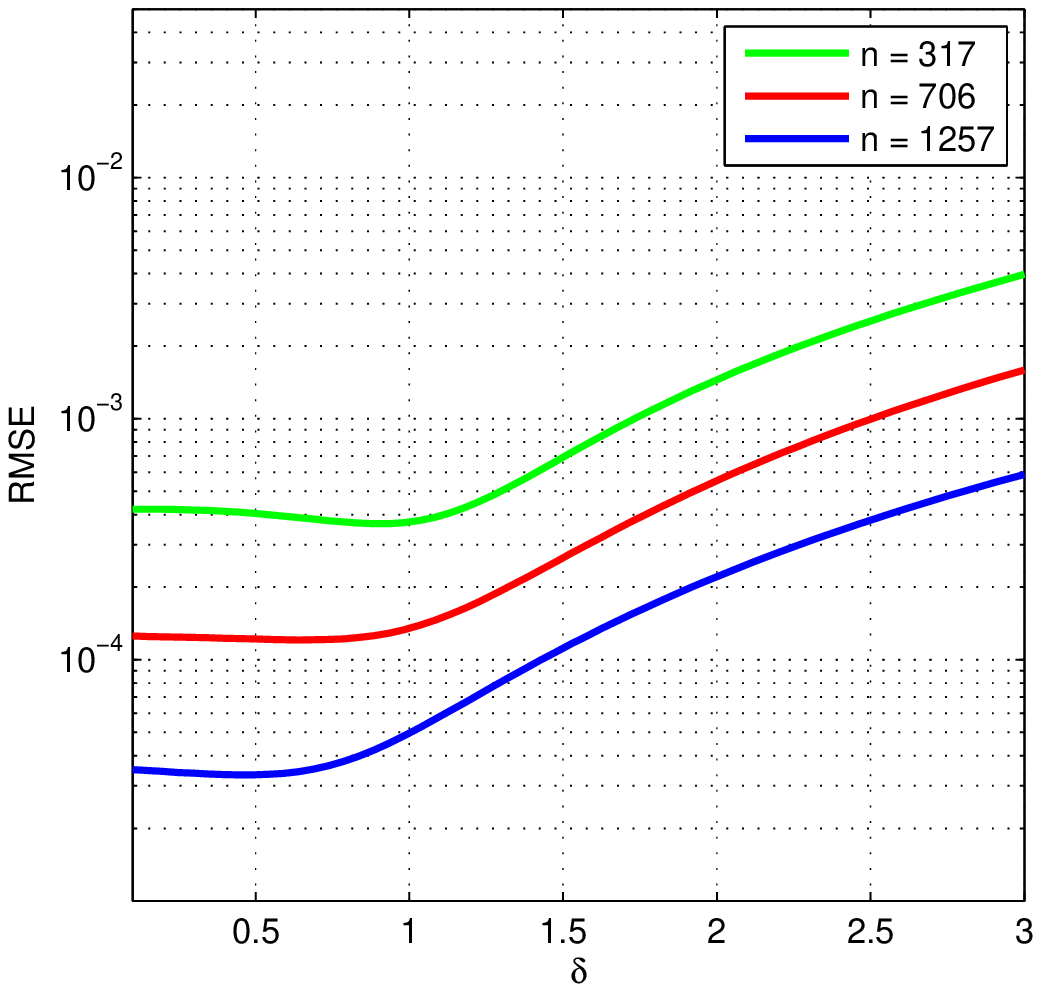} \hskip -0.85cm 
\includegraphics[height=3.4cm]{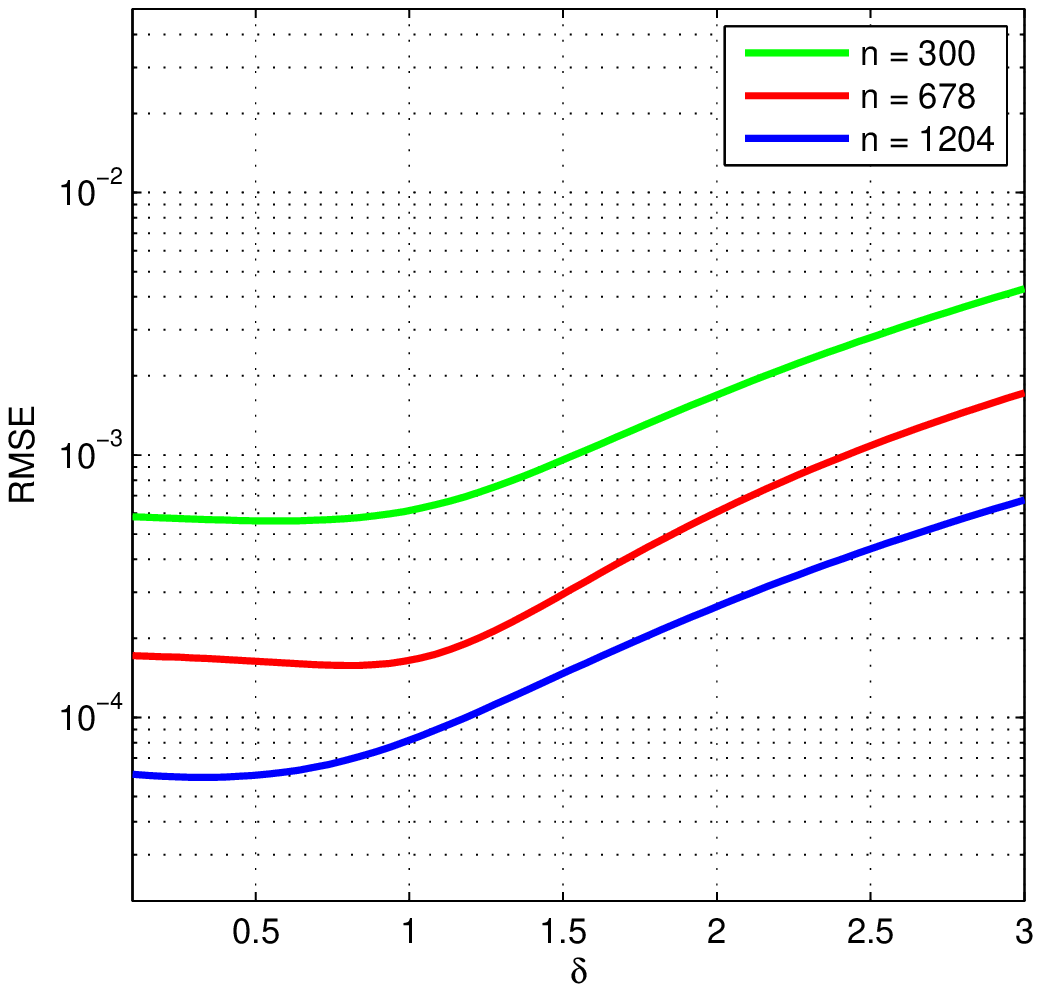}
\caption{RMSEs by varying the value of $\delta$. Left: triangle; center: disk; right: hexagon.}
\label{fig_1_RMSE}
\end{figure}

After showing in Figure \ref{fig_1_RMSE} the stable behavior of RMSEs by varying the value of $\delta \in [0.1,3]$, for each of convex domains we report MAEs and RMSEs taking $\delta = 0.1$ as shape parameter of the Wendland $C^2$ function. Then, since we are also concerned to point out the efficiency of the algorithm, in Tables \ref{tab_1}--\ref{tab_3} we show CPU times computed in seconds. 

\begin{table} 
\caption{Errors and CPU times (in seconds) for triangle using $\delta = 0.1$.}
\begin{center}
\begin{tabular}{cccc}
\hline\noalign{\smallskip}
$n$ & MAE & RMSE & time \\
\noalign{\smallskip}
\hline
\noalign{\smallskip}
$51$  & $1.04{\rm E}-01$ & $1.06{\rm E}-02$ & 0.1 \\ 
$200$  & $6.57{\rm E}-02$ & $3.60{\rm E}-03$ & 0.2 \\ 
$451$  & $1.26{\rm E}-02$ & $6.11{\rm E}-04$ & 0.3 \\ 
$805$  & $7.39{\rm E}-03$ & $3.00{\rm E}-04$ & 0.4 \\ 
$1256$  & $3.72{\rm E}-03$ & $1.65{\rm E}-04$ & 0.6 \\ 
\hline
\end{tabular}
\end{center} \label{tab_1}
\end{table}

\begin{table} 
\caption{Errors and CPU times (in seconds) for disk using $\delta = 0.1$.}
\begin{center}
\begin{tabular}{cccc}
\hline\noalign{\smallskip}
$n$ & MAE & RMSE & time \\
\noalign{\smallskip}
\hline
\noalign{\smallskip}
$80$  & $2.64{\rm E}-02$ & $4.94{\rm E}-03$ & 0.1 \\ 
$317$  & $5.12{\rm E}-03$ & $4.22{\rm E}-04$ & 0.2 \\ 
$706$  & $1.99{\rm E}-03$ & $1.25{\rm E}-04$ & 0.4 \\ 
$1257$  & $3.29{\rm E}-04$ & $3.51{\rm E}-05$ & 0.6 \\ 
$1960$  & $3.23{\rm E}-04$ & $2.39{\rm E}-05$ & 0.9 \\ 
\hline
\end{tabular}
\end{center} \label{tab_2}
\end{table}

\begin{table} 
\caption{Errors and CPU times (in seconds) for hexagon using $\delta = 0.1$.}
\begin{center}
\begin{tabular}{cccc}
\hline\noalign{\smallskip}
$n$ & MAE & RMSE & time \\
\noalign{\smallskip}
\hline
\noalign{\smallskip}
$76$  & $4.43{\rm E}-02$ & $6.35{\rm E}-03$ & 0.1 \\ 
$300$  & $5.56{\rm E}-03$ & $5.82{\rm E}-04$ & 0.2 \\ 
$678$  & $2.38{\rm E}-03$ & $1.72{\rm E}-04$ & 0.4 \\ 
$1204$  & $6.40{\rm E}-04$ & $6.07{\rm E}-05$ & 0.6 \\ 
$1877$  & $6.32{\rm E}-04$ & $3.98{\rm E}-05$ & 0.8 \\ 
\hline
\end{tabular}
\end{center} \label{tab_3}
\end{table}

\begin{figure}
\centering
\includegraphics[height=4.9cm]{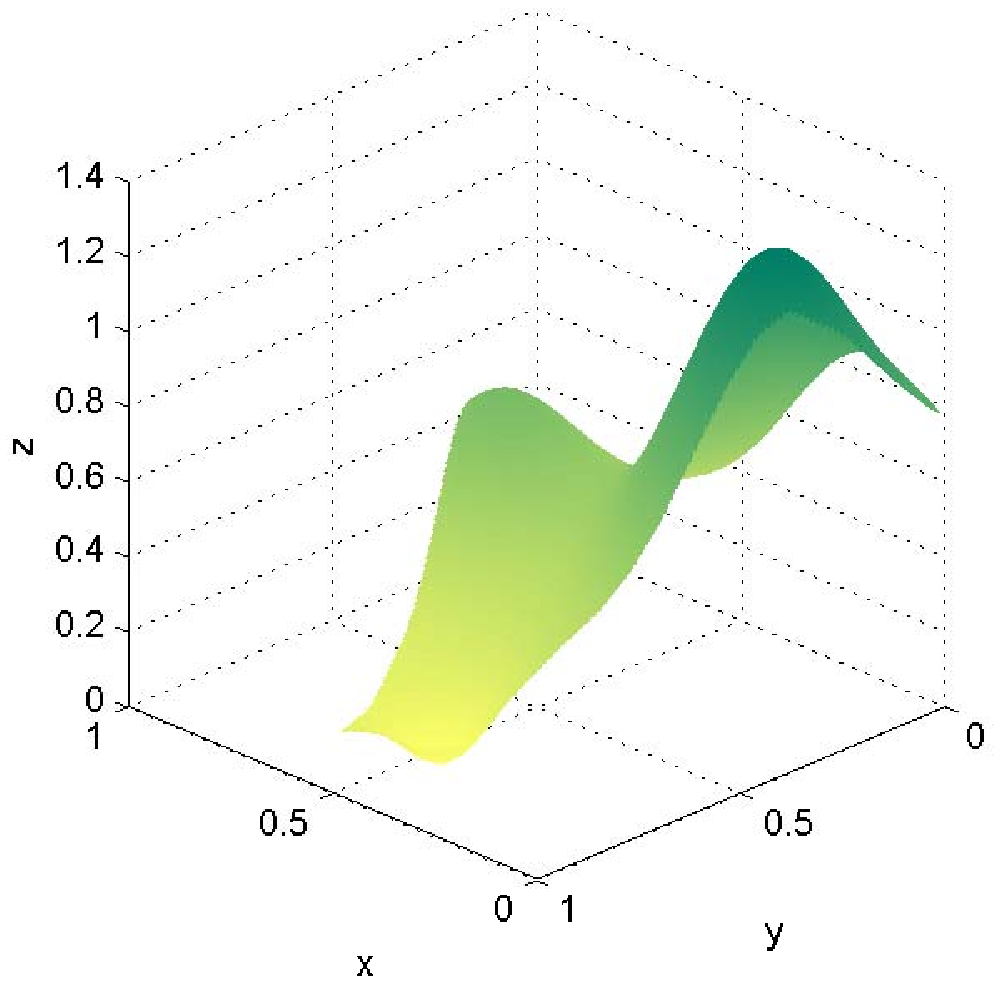} \hskip -1.4cm
\includegraphics[height=4.9cm]{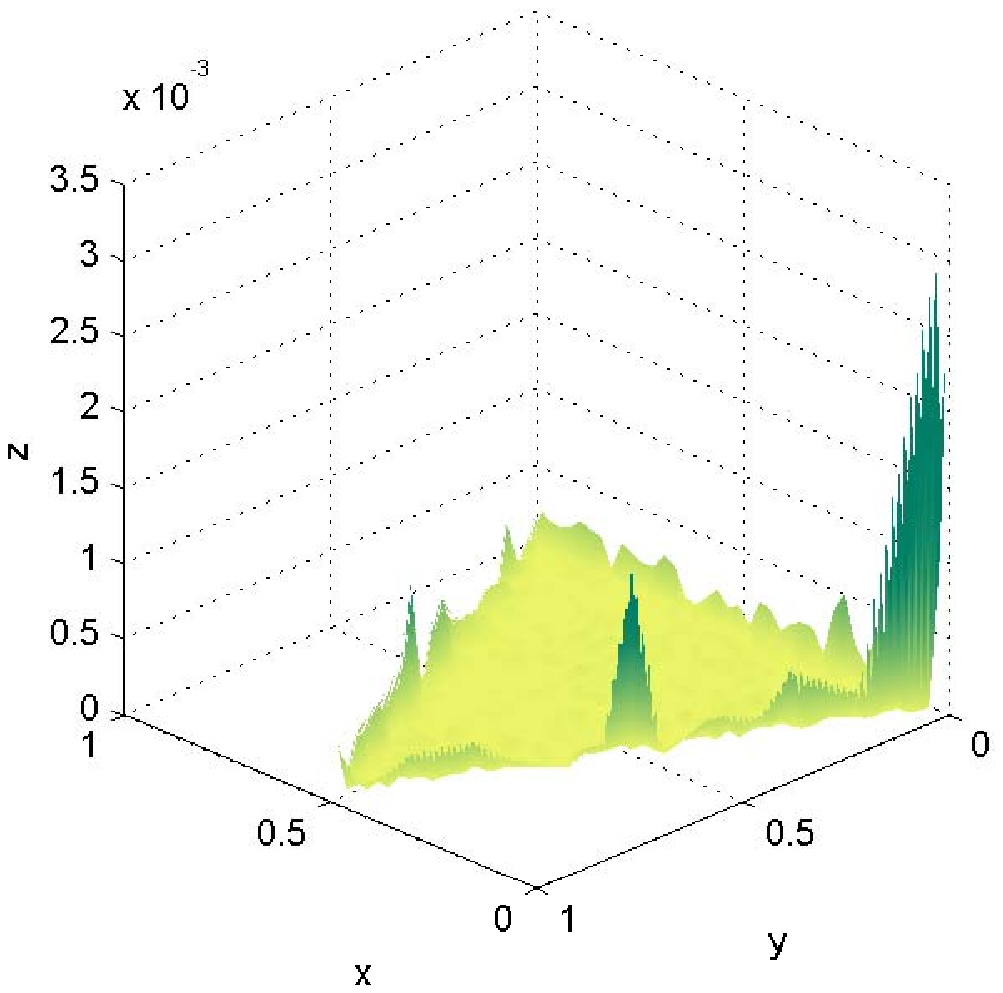}\\
\includegraphics[height=4.9cm]{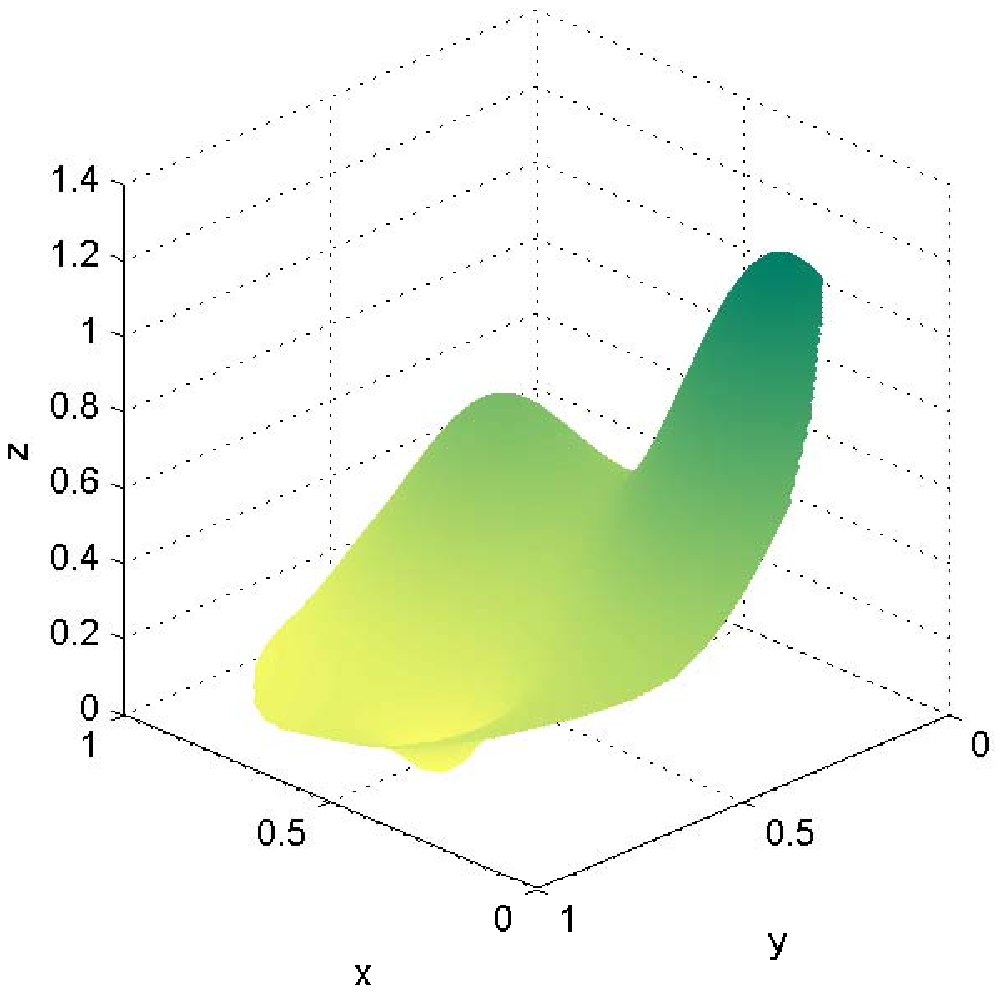} \hskip -1.4cm
\includegraphics[height=4.9cm]{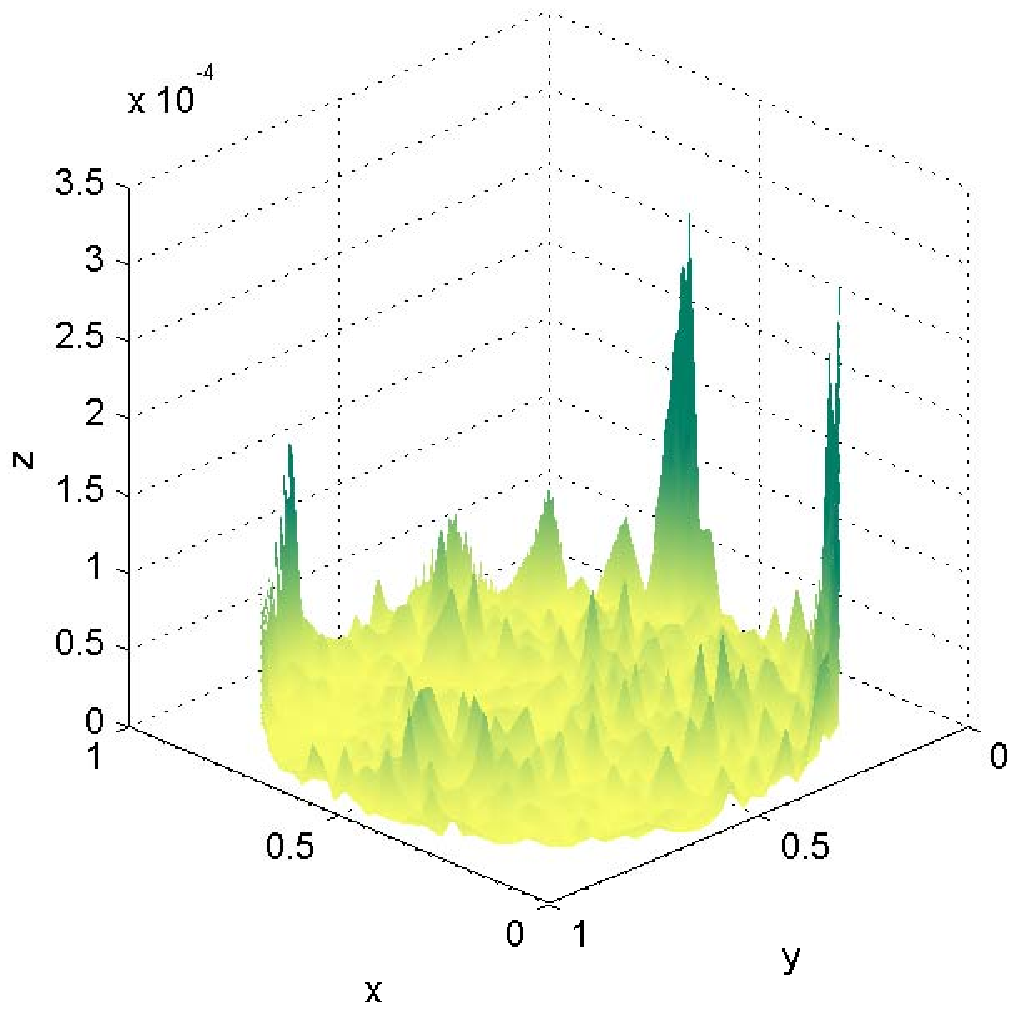}\\
\includegraphics[height=4.9cm]{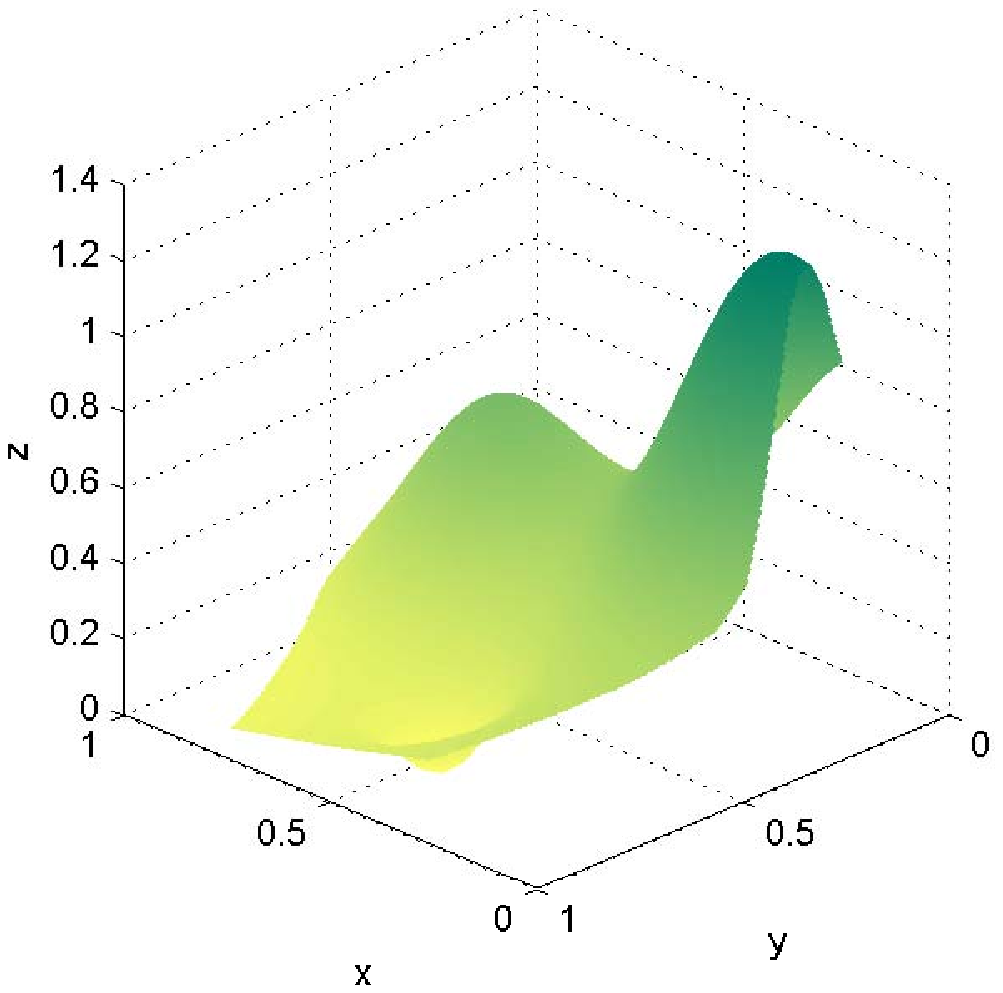} \hskip -1.4cm
\includegraphics[height=4.9cm]{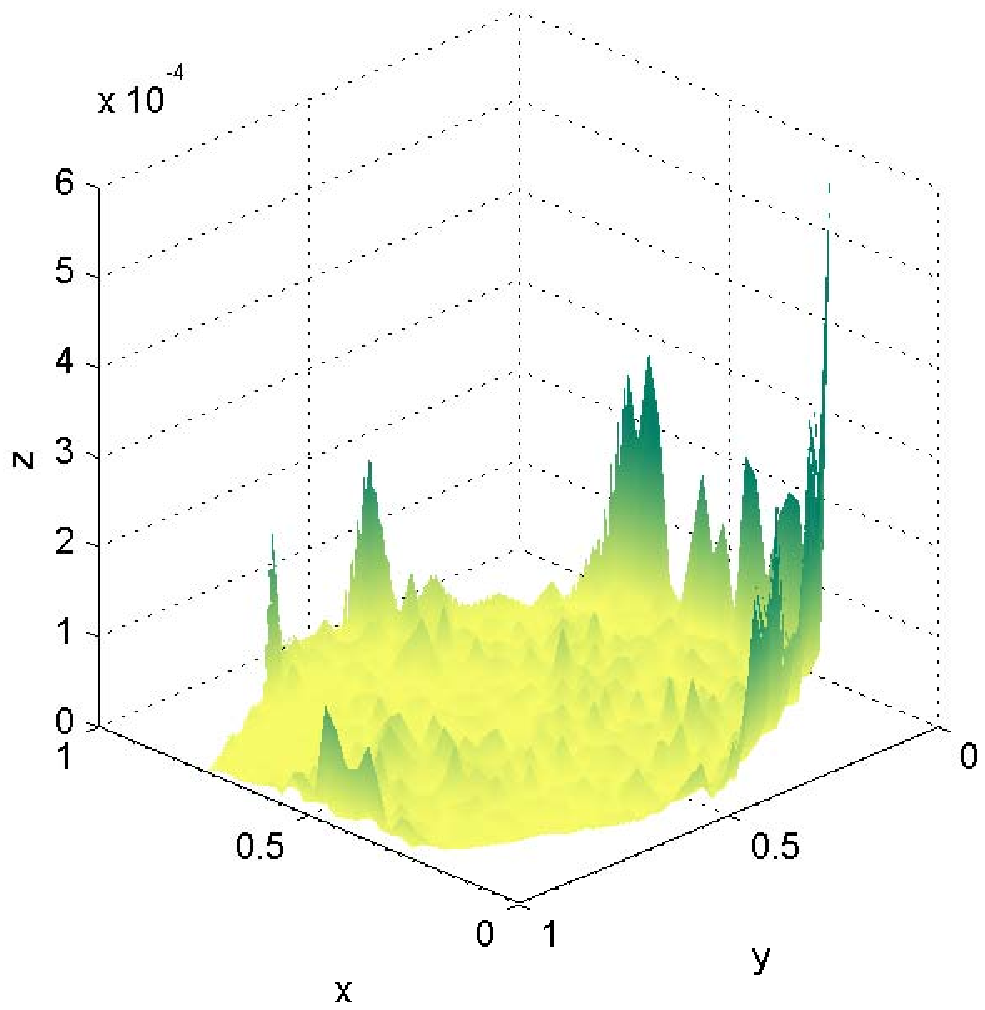}\\
\caption{2D Franke's function (left) and absolute errors (right) defined on convex domains. Top: triangle, $1256$ nodes; middle: disk, $1960$ nodes; bottom: hexagon, $1877$ nodes.}
\label{fig_1_franke}
\end{figure}

Finally, in Figure \ref{fig_1_franke} we represent the 2D Franke's function (left) and the absolute errors (right) computed on convex domains. This study shows that the maximum errors mainly concentrate on or close to the boundary of the convex hull. Note that, for shortness, in this paper we report numerical results obtained on a single example (or data set), but similar situations appear in all considered cases. 

\subsection{Results for 3D Convex Hulls} \label{num_exp_3D}

\begin{figure}
\centering
\includegraphics[height=3.5cm]{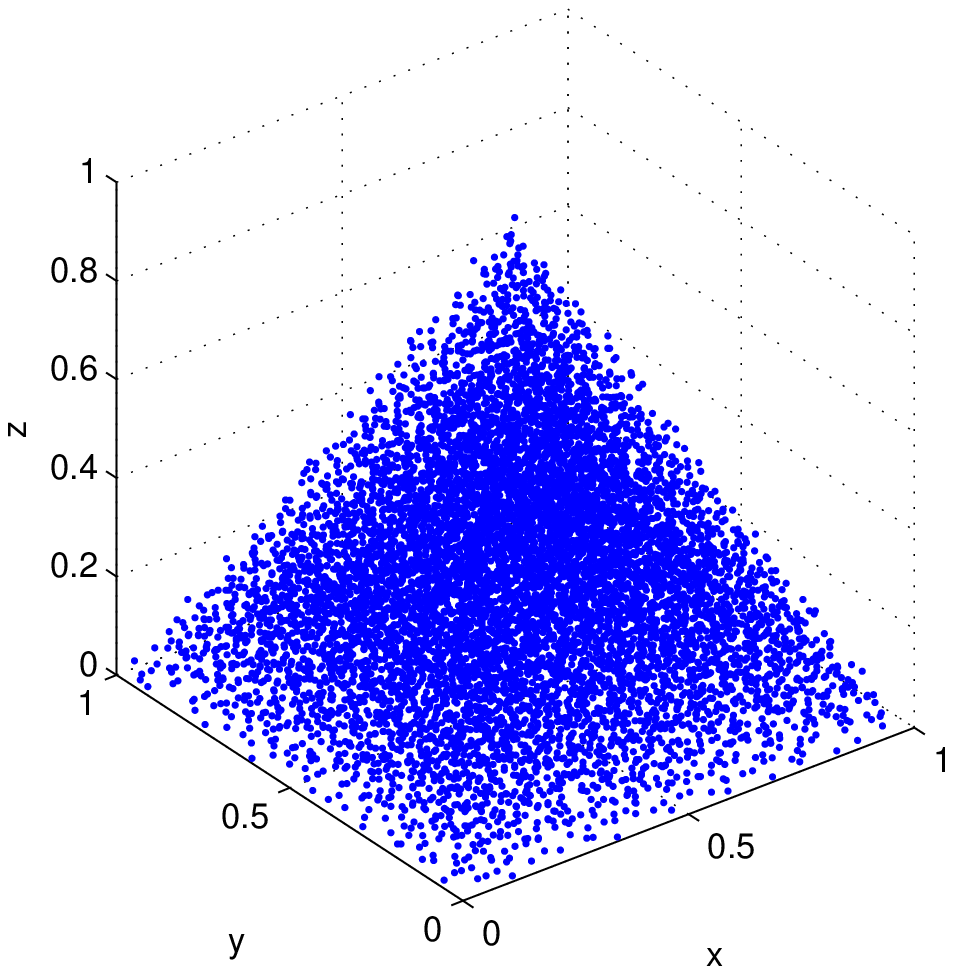} \hskip -1.cm
\includegraphics[height=3.5cm]{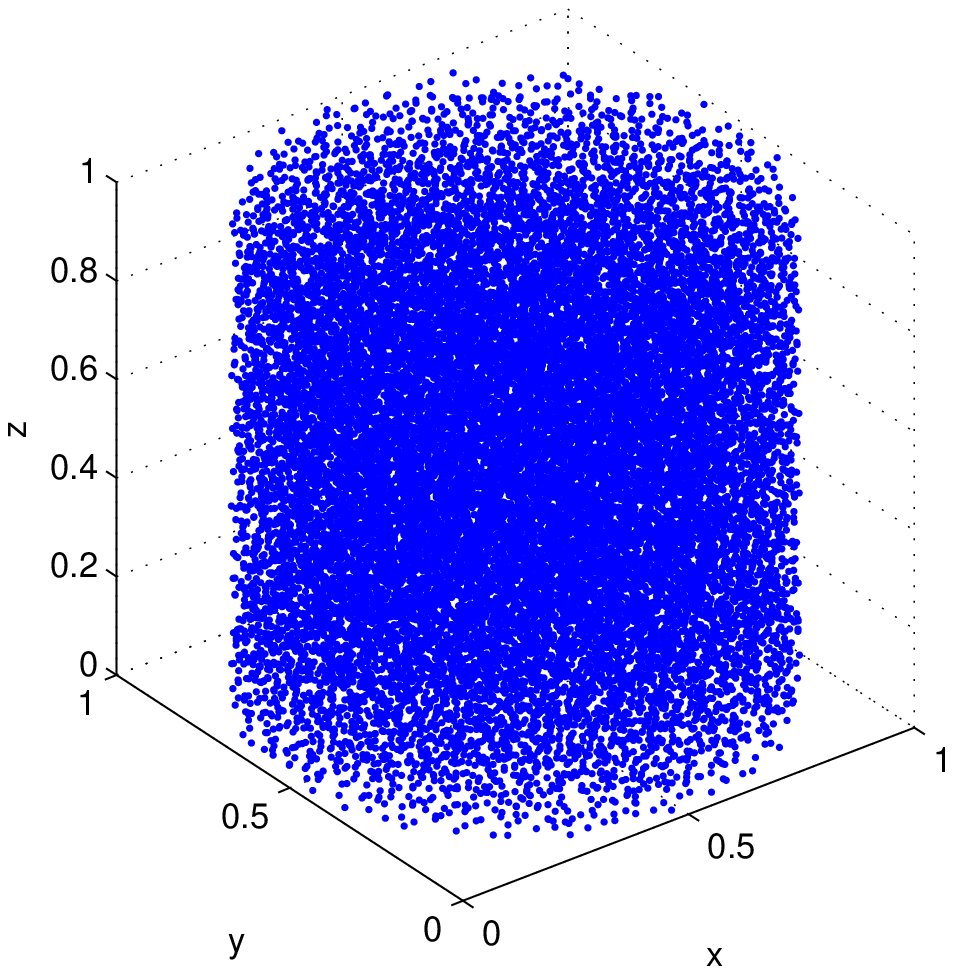} \hskip -1.cm 
\includegraphics[height=3.5cm]{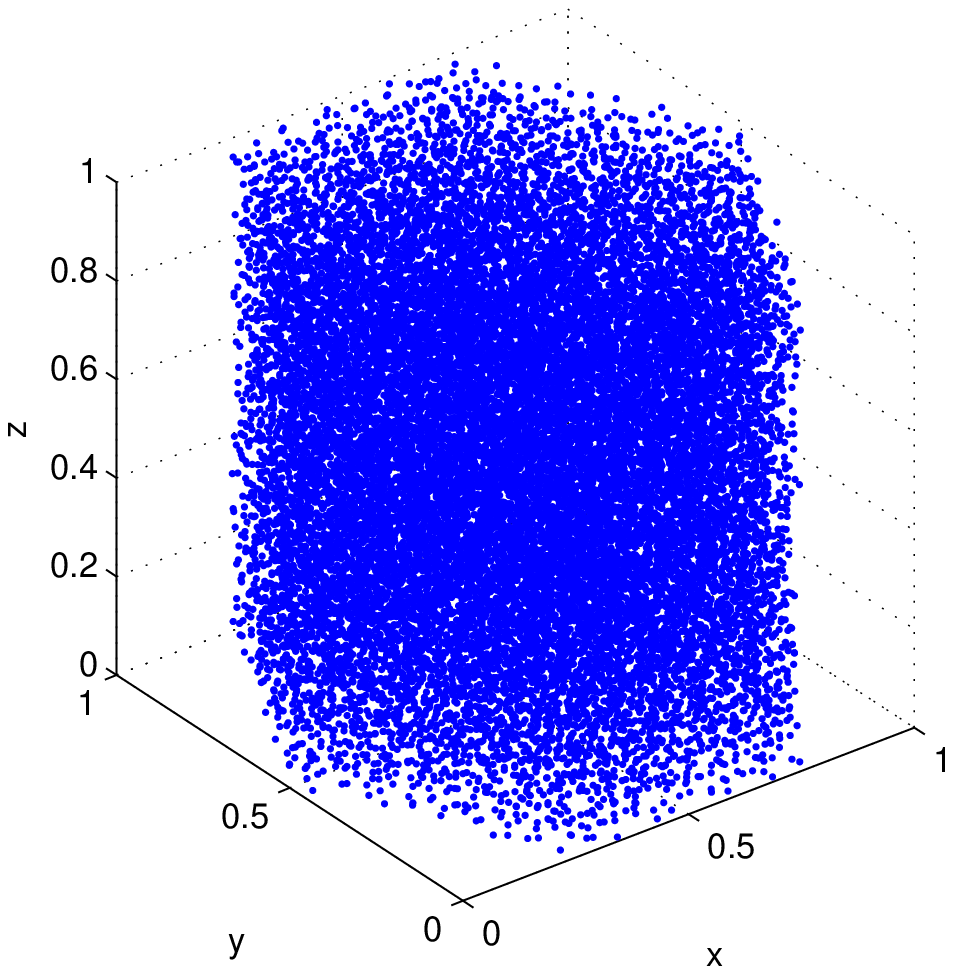}
\caption{Examples of points in 3D convex hulls. Left: pyramid, 8995 nodes; center: cylinder, 21177 nodes; right: hexagonal prism, 20249 nodes.}
\label{fig_2}
\end{figure}

In this subsection we instead report numerical results concerning trivariate interpolation. As earlier, we analyze accuracy and efficiency of the partition of unity algorithm for convex hulls, taking also in this case some sets of Halton scattered data points. Such points are located in three different convex domains: a pyramid, a cylinder and a hexagonal prism (see Figure \ref{fig_2}).

In the various tests we analyze the performances of the proposed algorithm taking the data values by 3D Franke's function, whose analytic expression is
\begin{eqnarray}
f_3(x_1,x_2,x_3)&=& \frac{3}{4}{\rm e}^{-\frac{(9x_1-2)^2+(9x_2-2)^2+(9x_3-2)^2}{4}}+\frac{3}{4} {\rm e}^{-\frac{(9x_1+1)^2}{49}-\frac{9x_2+1}{10}-\frac{9x_3+1}{10}} \nonumber \\
&+&\frac{1}{2} {\rm e}^{-\frac{(9x_1-7)^2+(9x_2-3)^2+(9x_3-5)^2}{4}}-\frac{1}{5} {\rm e}^{-(9x_1-4)^2-(9x_2-7)^2-(9x_3-5)^2}. \nonumber
\end{eqnarray}

As in the bivariate case, for each of convex hulls in Tables \ref{tab_4}--\ref{tab_6} we show MAEs, RMSEs and CPU times obtained by running our interpolation algorithm. These results are obtained taking $\delta = 0.1$. Here, we omit the graphs of RMSEs by varying $\delta$ because this study revealed a behavior similar to that outlined in Figure \ref{fig_1_RMSE}. Moreover, in dimension three we observed a even more relevant concentration of maximum errors on (or close to) the boundary of convex hulls.

\begin{table} 
\caption{Errors and CPU times (in seconds) for pyramid using $\delta = 0.1$.}
\begin{center}
\begin{tabular}{cccc}
\hline\noalign{\smallskip}
$n$ & MAE & RMSE & time \\
\noalign{\smallskip}
\hline
\noalign{\smallskip}
$335$  & $1.15{\rm E}-01$ & $5.03{\rm E}-03$ & 4.8 \\ 
$2670$  & $3.42{\rm E}-02$ & $5.88{\rm E}-04$ & 17.4 \\ 
$8995$  & $1.21{\rm E}-02$ & $1.60{\rm E}-04$ & 30.1 \\ 
$21337$  & $1.66{\rm E}-02$ & $1.59{\rm E}-04$ & 49.5 \\ 
$41665$  & $5.96{\rm E}-03$ & $5.95{\rm E}-05$ & 83.7 \\ 
\hline 
\end{tabular}
\end{center} \label{tab_4}
\end{table}

\begin{table} 
\caption{Errors and CPU times (in seconds) for cylinder using $\delta = 0.1$.}
\begin{center}
\begin{tabular}{cccc}
\hline\noalign{\smallskip}
$n$ & MAE & RMSE & time \\
\noalign{\smallskip}
\hline
\noalign{\smallskip}
$787$  & $3.47{\rm E}-01$ & $7.93{\rm E}-03$ & 15.2 \\ 
$6271$  & $9.91{\rm E}-03$ & $1.65{\rm E}-04$ & 44.6 \\ 
$21177$  & $2.00{\rm E}-03$ & $3.35{\rm E}-05$ & 77.6 \\ 
$50184$  & $1.48{\rm E}-03$ & $1.86{\rm E}-05$ & 130.7 \\ 
$97997$  & $8.58{\rm E}-04$ & $1.08{\rm E}-05$ & 209.0 \\ 
\hline
\end{tabular}
\end{center} \label{tab_5}
\end{table}

\begin{table} 
\caption{Errors and CPU times (in seconds) for hexagonal prism using $\delta = 0.1$.}
\begin{center}
\begin{tabular}{cccc}
\hline\noalign{\smallskip}
$n$ & MAE & RMSE & time \\
\noalign{\smallskip}
\hline
\noalign{\smallskip}
$754$  & $1.65{\rm E}-01$ & $4.49{\rm E}-03$ & 15.1 \\ 
$6002$  & $7.25{\rm E}-03$ & $1.33{\rm E}-04$ & 44.5 \\ 
$20249$  & $3.02{\rm E}-03$ & $3.78{\rm E}-05$ & 77.3 \\ 
$47997$  & $1.78{\rm E}-03$ & $2.00{\rm E}-05$ & 127.0 \\ 
$93754$  & $1.15{\rm E}-03$ & $1.27{\rm E}-05$ & 202.9 \\ 
\hline
\end{tabular}
\end{center} \label{tab_6}
\end{table}


\section{Conclusions and Future Work} \label{concl}

In this paper we presented a new algorithm for multivariate interpolation of scattered data sets lying in convex domains (or hulls) $\Omega \subseteq \RR^N$, for any $N \geq 2$. It is based on the partition of unity interpolation using local RBF interpolants and compactly supported weight functions. To partition the points in $\Omega$, we used a kd-tree data structure efficiently answering the range search query. 

As future work, we expect to build new data structures for partitioning data in convex hulls using efficient cell-based searching procedures. The new code should allow us to further reduce CPU times, making it suitable and applicable in several fields of applied mathematics and scientific computing.



\subsubsection*{Acknowledgments.} The first author acknowledges financial support from the GNCS--INdAM and the University of Torino via grant \lq\lq Approssimazione di dati sparsi e sue applicazioni\rq\rq. 

\end{document}